\documentclass[12pt]{article}
\usepackage{amssymb,amsmath}

\newcommand{\Z}{{\mathbb Z}}

\newcommand{\Q}{{\mathbb Q}}

\newcommand{\Fps}{{\mathbb F}_{p^2}}
\newcommand{\Zps}{{\mathbb Z}_{p^2}}
\newcommand{\Qps}{{\mathbb Q}_{p^2}}
\newcommand{\Fq}{{\mathbb F}_q}

\newcommand{\OO}{{\cal O}}
\newcommand{\M}{\mathcal{M}}

\newcommand{\Gal}{\text{Gal}}
\newcommand{\Res}{\text{Res}}
\newcommand{\id}{\text{id}}
\newcommand{\N}{\text{N}}
\newcommand{\Tr}{\text{Tr}}
\newcommand{\Rbar}{\overline{R}}

\newcommand{\proof}{\noindent{\em Proof: }}
\newcommand{\qed}{\hspace{\fill}$\square$}
\newcommand{\ra}{\rightarrow}

\newcommand{\tst}{\textstyle}
\newcommand{\smallfrac}[2]{{\textstyle\frac{#1}{#2}}}
\newcommand{\mubold}{{\boldsymbol\mu}}

\newtheorem{theorem}{Theorem}
\newtheorem{lemma}[theorem]{Lemma}
\newtheorem{prop}[theorem]{Proposition}
\newtheorem{cor}[theorem]{Corollary}
\newenvironment{remark}{\noindent\refstepcounter{theorem}{\bf
Remark \thesection.\arabic{theorem}} }{}

\numberwithin{equation}{section}
\numberwithin{theorem}{section}

\title{Indices of inseparability and refined ramification
breaks}
\author{Kevin Keating \\
Department of Mathematics \\
University of Florida \\
Gainesville, FL 32611 \\
USA \\[.2cm]
{\tt keating@ufl.edu}}

\begin{document}

\maketitle

\begin{abstract}
\noindent
Let $K$ be a finite extension of $\Q_p$ which contains a
primitive $p$th root of unity $\zeta_p$.  Let $L/K$ be
a totally ramified $(\Z/p\Z)^2$-extension which has a
single ramification break $b$.  In \cite{necessity}
Byott and Elder defined a ``refined ramification break''
$b_*$ for $L/K$.  In this paper we prove that if $p>2$
and the index of inseparability $i_1$ of $L/K$ is not
equal to $p^2b-pb$ then ${b_*=i_1-p^2b+pb+b}$.
\end{abstract}

\section{Introduction}

Let $K$ be a finite extension of $\Q_p$, let $L/K$ be a
finite Galois extension, and let $\pi_L$ be a
uniformizer for $L$.  For simplicity we assume that
$L/K$ is a totally ramified extension of degree $p^n$
for some $n\ge1$.  The (lower) ramification breaks of
$L/K$ are the integers $v_L(\sigma(\pi_L)-\pi_L)-1$ for
$\sigma\in\Gal(L/K)$, $\sigma\not=\id_L$.  The extension
$L/K$ has at most $n$ distinct ramification breaks; if
there are fewer than $n$ breaks then $L/K$ may be viewed
as having degenerate ramification data.

     There have been several attempts to supply the
``missing'' ramification data in the cases where $L/K$
has fewer than $n$ breaks.  The indices of
inseparability $i_0,i_1,\dots,i_n$ of $L/K$ were defined
by Fried \cite{fried} in characteristic $p$ and by
Heiermann \cite{heier} in characteristic 0.  The indices
of inseparability determine the ramification breaks of
$L/K$ in all cases.  As for the opposite direction, if
$L/K$ has $n$ distinct ramification
breaks then the breaks determine the indices of
inseparability, but if $L/K$ has fewer than
$n$ breaks then the indices of inseparability are not
completely determined by the breaks.  Thus the indices
of inseparability give extra information about the
extension $L/K$ which can be viewed as the missing
ramification data.

     In \cite{new,necessity}, Byott and Elder
described an alternative method for supplying missing
ramification data by
defining refined lower ramification breaks for extensions
with fewer than $n$ ordinary breaks.  Suppose $L/K$ is a
totally ramified $(\Z/p\Z)^2$-extension with a single
(ordinary) ramification break $b$.  Then $L/K$ has one
refined break $b_*$, which is computed in
\cite{necessity} under the assumption that $K$ contains
a primitive $p$th root of unity.  Byott and Elder also
showed that the Galois module
structure of $\OO_L$ determines $b_*$ in certain cases.

     In this paper we study the relationship between the
index of inseparability $i_1$ of $L/K$ and the refined
ramification break $b_*$.  In particular, when $p>2$ and
$i_1\not=p^2b-pb$ we
give a formula which expresses $b_*$ in terms of $i_1$.
Our approach is based on the methods given
in \cite{elem} for computing
$i_1$ in terms of the norm group
$\N_{L/K}(L^{\times})$.
We relate these methods to the Byott-Elder formula
for $b_*$ using Vostokov's formula \cite{vos} for
computing the Kummer pairing $\langle\;,\,\rangle_p:
K^{\times}\times K^{\times}\ra\mubold_p$.  The
calculations are simplified somewhat through the use
of the Artin-Hasse exponential series $E_p(X)$.

     The author would like to thank the referee for
writing a very careful and thorough review of this paper.
\\[\medskipamount]
{\bf \large Notation}
\\[\smallskipamount]
$K=$ finite extension of $\Q_p$. \\
$K_0/\Q_p=$ maximum unramified subextension of $K/\Q_p$. \\
$v_K=$ valuation on $K$ normalized so that
$v_K(K^{\times})=\Z$. \\
$e=v_K(p)=$ absolute ramification index of $K$. \\
$\OO_K=\{\alpha\in K:v_K(\alpha)\ge0\}=$ ring of integers
of $K$. \\
$\M_K=\{\alpha\in K:v_K(\alpha)\ge1\}=$ maximal ideal
of $\OO_K$. \\
$\Fq\cong\OO_K/\M_K=$ residue field of $K$. \\
$U_K^c=1+\M_K^c$ for $c\ge1$. \\
$K^{ab}=$ maximal abelian extension of $K$. \\
$L/K=$ totally ramified $(\Z/p\Z)^2$-subextension of
$K^{ab}/K$ with one ramification break $b$. \\
$\pi_L=$ uniformizer for $L$. \\
$\pi_K=\N_{L/K}(\pi_L)=$ uniformizer for $K$. \\
$\zeta_n=$ primitive $n$th root of unity in $K^{ab}$. \\
$\mubold_n=\langle\zeta_n\rangle$. \\
$\Zps=\Z_p[\mubold_{p^2-1}]$.

\section{The Artin-Hasse exponential series and
truncated exponentiation}

     In this section we study the relation between the
Artin-Hasse exponential series and the ``truncated
exponentiation'' polynomials of Byott-Elder.  We also
use the Artin-Hasse exponential series to obtain a new
version of a formula from \cite{elem} for the index of
inseparability $i_1$ of a $(\Z/p\Z)^2$-extension with a
single ramification break.

     The Artin-Hasse exponential series is defined by 
\begin{equation} \label{Ep}
E_p(X)=\exp\left(X+\frac{1}{p}X^p+\frac{1}{p^2}X^{p^2}
+\cdots\right),
\end{equation}
where $\exp(X)\in\Q[[X]]$ is the usual exponential
series.  Let $\mu$ denote the M\"obius function.  Then
by Lemma 9.1 in \cite[I]{FV} we have
\begin{equation} \nonumber
E_p(X)=\prod_{p\nmid c}\,(1-X^c)^{-\mu(c)/c}.
\end{equation}
Thus the coefficients of $E_p(X)$ lie in
$\Z_{(p)}=\Q\cap\Z_p$.
For each $i\ge1$ the power series $E_p(X)=1+X+\cdots$
induces a bijection from $\M_K^i$ onto $U_K^i$.  For
$\kappa,\lambda\in\M_K$ we have $E_p(\kappa)\equiv
E_p(\lambda)\pmod{\M_K^i}$ if and only if
$\kappa\equiv\lambda\pmod{\M_K^i}$.   Let
$\Lambda_p:U_K^1\ra\M_K$ denote the inverse of the
bijection from $\M_K$ to $U_K^1$ induced by $E_p(X)$.
Then for $u,v\in U_K^1$ we have
$\Lambda_p(u)\equiv\Lambda_p(v)\pmod{\M_K^i}$
if and only if $u\equiv v\pmod{\M_K^i}$.

     Let $\psi(X)\in XK[[X]]$ and $\alpha\in K$.  The
$\alpha$ power of $1+\psi(X)$ is a series in $K[[X]]$
defined by
\begin{equation} \nonumber
(1+\psi(X))^{\alpha}=\sum_{n=0}^{\infty}\,
\binom{\alpha}{n}\psi(X)^n,
\end{equation}
where
\begin{equation} \nonumber
\binom{\alpha}{n}=\frac{\alpha(\alpha-1)(\alpha-2)\dots
(\alpha-(n-1))}{n!}.
\end{equation}
Motivated by this formula, Byott and Elder
\cite[1.1]{new} defined truncated exponentiation by
\begin{equation} \nonumber
(1+\psi(X))^{[\alpha]}=\sum_{n=0}^{p-1}\,
\binom{\alpha}{n}\psi(X)^n.
\end{equation}
Thus $(1+X)^{[\alpha]}$ is a polynomial with
coefficients in $K$; if $\alpha\in\OO_K$ then the
coefficients of $(1+X)^{[\alpha]}$ lie in $\OO_K$.  For
$u\in U_K^1$ define $u^{[\alpha]}$ to be the value of
$(1+X)^{[\alpha]}$ at $X=u-1$.

\begin{lemma} \label{alpha}
Let $\alpha\in K$.  Then
$E_p(X)^{[\alpha]}\equiv E_p(\alpha X)\pmod{X^p}$.
\end{lemma}

\proof We have $E_p(X)^{[\alpha]}\equiv\exp(X)^{\alpha}
\equiv\exp(\alpha X)\equiv E_p(\alpha X)\pmod{X^p}$. \qed

\begin{prop} \label{hom}
Let $i\ge1$, let $u,v\in U_K^i$, and let
$\alpha\in\OO_K$.  Then 
\begin{alignat*}{2}
\Lambda_p(uv)&\equiv\Lambda_p(u)+\Lambda_p(v)
&&\pmod{\M_K^{pi}} \\
\Lambda_p(u^{[\alpha]})&\equiv\alpha\Lambda_p(u)
&&\pmod{\M_K^{pi}}.
\end{alignat*}
\end{prop}

\proof Set $\kappa=\Lambda_p(u)$ and $\lambda=\Lambda_p(v)$.
Then $\kappa,\lambda\in\M_K^i$, so by equation (6) in
\cite[p.\,52]{pdg} we have
\begin{equation} \nonumber
E_p(\kappa)E_p(\lambda)\equiv E_p(\kappa+\lambda)
\pmod{\M_K^{pi}}.
\end{equation}
In addition, by Lemma~\ref{alpha} we get
\begin{equation} \nonumber
E_p(\kappa)^{[\alpha]}\equiv E_p(\alpha\kappa)
\pmod{\M_K^{pi}}.
\end{equation}
Applying $\Lambda_p$ to these congruences gives the
desired results. \qed

\begin{cor} \label{OK}
Let $i\ge1$.  The scalar multiplication
$\alpha\cdot u=u^{[\alpha]}$ induces an $\OO_K$-module
structure on the group $U_K^i/U_K^{pi}$.  Furthermore,
$\Lambda_p$ induces an isomorphism of $\OO_K$-modules
from $U_K^i/U_K^{pi}$ onto $\M_K^i/\M_K^{pi}$.
\end{cor}

\begin{cor}
Let $u\in U_K^i$ and $\alpha\in\Z_p$.  Then
$u^{\alpha}\equiv u^{[\alpha]}\pmod{\M_K^{pi}}$.
\end{cor}

\proof For $n\ge1$ we have $\Lambda_p(u^n)\equiv
n\Lambda_p(u)\equiv\Lambda_p(u^{[n]})\pmod{\M_K^{pi}}$.
\qed

\begin{cor} \label{Zp}
Let $i\ge1$ and let $A$ be a subgroup of $U_K^i$ which
contains $U_K^{pi}$.  Then $\Lambda_p(A)$ is a
$\Z_p$-module.
\end{cor}

\begin{cor} \label{sub}
Let $i\ge1$ and let $A,B$ be subgroups of $U_K^i$ such
that $U_K^{pi}\subset B$.  Then
$\Lambda_p(AB)=\Lambda_p(A)+\Lambda_p(B)$.
\end{cor}

\proof We clearly have
$\Lambda_p(AB)\supset\Lambda_p(A)$ and
$\Lambda_p(AB)\supset\Lambda_p(B)$.  Hence by
Corollary~\ref{Zp} we get
$\Lambda_p(AB)\supset\Lambda_p(A)+\Lambda_p(B)$.  Let
$a\in A$, $b\in B$.  Then
$\Lambda_p(ab)=\Lambda_p(a)+\Lambda_p(b)+m$ for some
$m\in\M_K^{pi}$.  Let $b'\in U_K^i$ be such that
$\Lambda_p(b')=\Lambda_p(b)+m$.  Then $b\equiv
b'\pmod{\M_K^{pi}}$, so $b'\in B$.  Hence
$\Lambda_p(AB)\subset\Lambda_p(A)+\Lambda_p(B)$.  We
conclude that
$\Lambda_p(AB)=\Lambda_p(A)+\Lambda_p(B)$.~\qed \medskip

     Let $\Qps=\Q_p(\zeta_{p^2-1})$ denote the
unramified extension of $\Q_p$ of degree 2, and let
$\Zps$ denote the ring of integers of $\Qps$.

\begin{cor} \label{module}
Assume $\mubold_{p^2-1}\subset K$ and let $A$ be a
subgroup of $U_K^i$ which contains $U_K^{pi}$.  Then
$\Lambda_p(A)$ is a $\Zps$-module if and only if $A$ is
stable under the map $a\mapsto a^{[\eta]}$ for every
$\eta\in{\mubold}_{p^2-1}$.
\end{cor}

\proof This follows from Proposition~\ref{hom} and the
fact that $\Zps=\Z_p[\mubold_{p^2-1}]$. \qed

\begin{prop} \label{OK0}
Let $i,j$ be positive integers such that $pj\ge i$ and
$e+\lceil\frac{j}{p}\rceil\ge i$, and let $K_0/\Q_p$ be
the maximum unramified subextension of $K/\Q_p$.  Then
$\Lambda_p((K^{\times})^p\cap U_K^j)+\M_K^i$ is an
$\OO_{K_0}$-module.
\end{prop}

\proof If $i\le j$ then the claim is obvious, so we
assume $i\ge j+1$.  Then
\begin{equation} \nonumber
i\le e+\left\lceil\frac{i-1}{p}\right\rceil
\le e+\frac{i+p-2}{p}.
\end{equation}
It follows that $i\le\frac{pe}{p-1}+\frac{p-2}{p-1}$,
and hence that $i\le\lceil\frac{pe}{p-1}\rceil$.  By
applying Corollary~\ref{sub} with $i$ replaced by $j$,
$A=(K^{\times})^p\cap U_K^j$, and $B=U_K^i$ we get
\begin{equation} \nonumber
\Lambda_p(((K^{\times})^p\cap U_K^j)\cdot U_K^i)
=\Lambda_p((K^{\times})^p\cap U_K^j)+\M_K^i.
\end{equation}
Hence by Corollary~\ref{Zp} we see that
$\Lambda_p((K^{\times})^p\cap U_K^j)+\M_K^i$ is a
$\Z_p$-module.  Let $u\in(K^{\times})^p\cap U_K^j$ with
$c=v_K(u-1)<i$.  Then there is $\gamma\in\M_K$ such that
$u=E_p(\gamma)^p$.  Using (\ref{Ep}) we get
\begin{alignat*}{2}
u&\tst=\exp(p\gamma+\gamma^p+\frac1p\gamma^{p^2}+\dots) \\
&=\exp(p\gamma)\cdot E_p(\gamma^p).
\end{alignat*}
Since $c<\lceil\frac{pe}{p-1}\rceil$ and $c$ is an
integer we have $c<\frac{pe}{p-1}$, so $p\mid c$
and $v_K(\gamma)=\frac{c}{p}$.  Therefore $v_K(p\gamma)
=e+\frac{c}{p}\ge e+\lceil\frac{j}{p}\rceil\ge i$, and
hence $u\equiv E_p(\gamma^p)\pmod{\M_K^i}$.
On the other hand, for each $\gamma\in\M_K$ such that
$v_K(\gamma^p)\ge j$, the computations above show that
$E_p(\gamma^p)=E_p(\gamma)^p\cdot\exp(-p\gamma)$
lies in $((K^{\times})^p\cap U_K^j)\cdot U_K^i$.  It
follows that
\begin{equation} \label{mup}
\Lambda_p((K^{\times})^p\cap U_K^j)+\M_K^i
=\{\gamma^p:\gamma\in\M_K,\;v_K(\gamma^p)\ge j\}+\M_K^i.
\end{equation}
Let $q$ be the cardinality of the residue field of
$K$.  Then $\mubold_{q-1}\subset\OO_K$, so the right
side of (\ref{mup}) is stable under multiplication by
elements of $\mubold_{q-1}$.  Since
$\OO_{K_0}=\Z_p[\mubold_{q-1}]$, the proposition
follows.~\qed \medskip

\section{Two invariants of $L/K$}

Let $L/K$ be a totally ramified $(\Z/p\Z)^2$-extension
with a single ramification break $b$.  Then $1\le
b<\frac{pe}{p-1}$ and $p\nmid b$ (see for instance
\cite[p.\,398]{dalfurther}).  In this section we define
two further invariants of $L/K$: the refined
ramification break $b_*$ and the index of inseparability
$i_1$.  We also show how $i_1$ can be computed in terms
of the valuations of the coefficients of the minimum
polynomial over $K$ of a uniformizer for $L$.

     To motivate the definition of $b_*$ we first
reformulate the definition of $i(\sigma)$ for
$\sigma\in\Gal(L/K)$.  It is easily seen that
\begin{equation} \nonumber
i(\sigma)=\min\{v_L(\sigma(\alpha)-\alpha)-v_L(\alpha):
\alpha\in\OO_L,\;\alpha\not=0\}.
\end{equation}
Thus $i(\sigma)$ may be viewed as the valuation of
the operator $\sigma-1$ on $\OO_L$.  Now let
$\sigma_1,\sigma_2$ be generators for
$\Gal(L/K)\cong(\Z/p\Z)^2$.  Since $b$
is the unique ramification break of $L/K$, for $i=1,2$
we have $\sigma_i(\pi_L)-\pi_L=\beta_i$ with
$v_L(\beta_i)=b+1$.  Let $\delta\in\mubold_{q-1}$ be
such that $\beta_1/\beta_2\equiv\delta\pmod{\M_L}$.
Then
\begin{equation} \nonumber
\sigma_2^{[-\delta]}=\sum_{n=0}^{p-1}\,
\binom{-\delta}{n}(\sigma_2-1)^n
\end{equation}
is an element of the group ring $\OO_{K_0}[\Gal(L/K)]$.
We define
\begin{align*}
b_*&=\min\{v_L(\sigma_1\circ\sigma_2^{[-\delta]}(\alpha)
-\alpha)-v_L(\alpha):\alpha\in\OO_L,\;\alpha\not=0\}.
\end{align*}
Thus $b_*=i(\sigma_1\circ\sigma_2^{[-\delta]})$ is the
valuation of the operator
$\sigma_1\circ\sigma_2^{[-\delta]}-1$ on $\OO_L$.
It is proved in \cite{necessity} that $b_*$ does not
depend on the choice of generators $\sigma_1,\sigma_2$
for $\Gal(L/K)$.

     We now define the indices of inseparability of
$L/K$, following Heiermann \cite{heier}.  Let $\pi_L$ be
a uniformizer for $L$.  Then $\pi_K=\N_{L/K}(\pi_L)$
is a uniformizer for $K$, and there are unique
$c_h\in\mubold_{q-1}\cup\{0\}$ such that
\begin{equation} \nonumber
\pi_K=\sum_{h=0}^{\infty}\,c_h\pi_L^{h+p^2}.
\end{equation}
For $0\le j\le2$ set
\begin{align*}
i_j^*&=\min\{h\ge0:c_h\not=0,\;v_p(h+p^2)\le j\} \\
i_j&=\min\{i_{j'}^*+p^2e\cdot(j'-j):j\le j'\le2\}.
\end{align*}
Then $i_j^*$ may depend on the choice of $\pi_L$, but
$i_j$ does not (see \cite[Th.\,7.1]{heier}).
Furthermore, we have $0=i_2<i_1\le i_0$.  The relation
between the indices of inseparability and the ordinary
ramification data of $L/K$ is given by
\cite[Cor.\,6.11]{heier}.  In particular, we have
$i_0=p^2b-b$.

     As in \cite{elem} we let
\begin{equation} \nonumber
g(X)=X^{p^2}+a_1X^{p^2-1}+\dots+a_{p^2-1}X+a_{p^2}
\end{equation}
be the minimum polynomial for $\pi_L$ over $K$.  
Then by \cite[(3.5)]{elem} we get
\begin{align*}
i_1&=\min\left(\{p^2v_K(a_i)-i:1\le i\le p^2-1\}
\cup\{i_2+p^2e\}\right) \\
&=\min\left(\{p^2v_K(a_{pi})-pi:1\le i\le p-1\}
\cup\{i_2+p^2e,\,i_0\}\right) \\
&=\min\left(\{p^2v_K(a_{pi})-pi:1\le i\le p-1\}
\cup\{p^2e,\,p^2b-b\}\right).
\end{align*}
For $j>p^2$ write $j=p^2u+i$ with $1\le i\le p^2$ and set
$a_j=\pi_K^ua_i$.  Then
$v_K(a_{pi+p^2c})=v_K(a_{pi})+c$, so for every $l\ge0$
we have
\begin{align} \label{i1}
i_1&=\min\left(\{p^2v_K(a_{pi})-pi:l<i\le l+p,\,p\nmid i\}
\cup\{p^2e,\,p^2b-b\}\right).
\end{align}

     Let $H=\N_{L/K}(L^{\times})$ be the subgroup of
$K^{\times}$ which is associated to the abelian
extension $L/K$ by class field theory.  Since $b$ is the
only ramification break of $L/K$ we have
$U_K^{b+1}\le H$ and
\begin{equation} \label{break}
U_K^b/(H\cap U_K^b)\cong K^{\times}/H\cong\Gal(L/K).
\end{equation}

\begin{theorem} \label{combining}
Let $p>2$, let $L/K$ be a totally ramified
$(\Z/p\Z)^2$-extension with a single ramification break
$b\ge1$, and set $H=\N_{L/K}(L^{\times})$.  If
$\mubold_{p^2-1}\not\subset K$ let $k=b$; otherwise let
$k$ be the smallest nonnegative integer such that
$\Lambda_p(H\cap U_K^{k+1})$ is a $\Zps$-module.  Then
\begin{equation} \nonumber
i_1=\min\{p^2b-pk,\,p^2e,\,p^2b-b\}.
\end{equation}
\end{theorem}

\proof Let $i\ge1$ satisfy $p\nmid i$.  Then by
\cite[(3.25)]{elem} we have
\begin{equation} \nonumber
\N_{L/K}(E_p(r\pi_L^i))\equiv
E_p(\pi_K^ir^{p^2})\cdot E_p(-ia_{pi}r^p-ia_ir)
\pmod{\M_K^{b+1}}.
\end{equation}
By \cite[Lemma 3.2]{elem} we have
\begin{align} \nonumber
v_K(a_i)&\ge b-\frac{b-i}{p^2}
=\left(1-\frac{1}{p^2}\right)b+\frac{1}{p^2}\cdot i \\
v_K(a_{pi})&\ge b-\frac{pb-pi}{p^2}
=\left(1-\frac{1}{p}\right)b+\frac{1}{p}\cdot i.
\label{api}
\end{align}
Hence if $i\le b$ then $v_K(a_i)\ge i$ and
$v_K(a_{pi})\ge i$, with strict inequalities if $i<b$.
It follows that
\begin{equation} \label{NEp}
\N_{L/K}(E_p(r\pi_L^i))\equiv E_p(\beta_i(r))
\pmod{\M_K^{b+1}},
\end{equation}
with $\beta_i(r)=\pi_K^ir^{p^2}-ia_{pi}r^p-ia_ir$.
In addition, we have $v_K(\beta_i(r))\ge i$, with
equality if $i<b$ and $r\not=0$.

     Since $\Lambda_p(H\cap U_K^{b+1})=\M_K^{b+1}$ we
have $k\le b$.  We claim that $v_K(a_{pi})\ge b+1$ for
all $i\ge k+1$ such that $p\nmid i$.  If $k=b$ this
follows from (\ref{api}).  Let $k<b$ and suppose the
claim is false.  Let $h\ge k+1$ be maximum with the
property that $p\nmid h$ and $v_K(a_{ph})\le b$.  Since
$a_{p(h+p)}=\pi_Ka_{ph}$ we see that a maximum $h$
exists, and that $v_K(a_{ph})=b$.  Since
$H\cap U_K^{k+1}\supset U_K^{b+1}$, it follows from
(\ref{NEp}) and Corollary~\ref{sub} that
$E_p(\beta_h(r))\in H\cap U_K^{k+1}$ for all
$r\in\mubold_{q-1}\cup\{0\}$.  By the definition of $k$,
$\Lambda_p(H\cap U_K^{k+1})$ is a $\Zps$-module.  Hence
for every $r\in{\mubold}_{q-1}$ and
$\eta\in{\mubold}_{p^2-1}$,
\begin{equation} \nonumber
\eta\beta_h(r)-\beta_h(\eta r)=ha_{ph}r^p(\eta^p-\eta)
\end{equation}
lies in $\Lambda_p(H\cap U_K^{k+1})$.  Since every coset
of $\M_K^{b+1}$ in $\M_K^b$ is represented by an element
of this form, and
\begin{equation} \nonumber
\Lambda_p(H\cap U_K^{k+1})\supset
\Lambda_p(U_K^{b+1})=\M_K^{b+1}, 
\end{equation}
it follows that
$\Lambda_p(H\cap U_K^{k+1})\supset\M_K^b$.  Hence
$H\supset E_p(\M_K^b)=U_K^b$, which contradicts
(\ref{break}).  This proves our claim, so we have
\begin{equation} \label{star}
p^2b-pk\le p^2v_K(a_{pi})-pi
\end{equation}
for all $i$ such that $k<i\le k+p$ and $p\nmid i$.

     Set $m=\min\{p^2b-pk,\,p^2e,\,p^2b-b\}$.  Suppose
$m=p^2b-b$.  Then $k\le\frac{b}{p}$, so by the preceding
paragraph we have $v_K(a_{pi})\ge b+1$ for all
$i>\frac{b}{p}$ such that $p\nmid i$.  Hence by
(\ref{i1}) we get
\begin{align*}
i_1&=\min(\{p^2v_K(a_{pi})-pi:\smallfrac{b}{p}<i\le
\smallfrac{b}{p}+p,\,p\nmid i\}
\cup\{p^2e,\,p^2b-b\}) \\
&=p^2b-b.
\end{align*}

     Suppose $m=p^2e$.  Then $k\le p(b-e)$, so
$v_K(a_{pi})\ge b+1$ for all $i>p(b-e)$ such that
$p\nmid i$.  Hence by (\ref{i1}) we have
\begin{align*}
i_1&=\min(\{p^2v_K(a_{pi})-pi:p(b-e)<i<p(b-e)+p\}
\cup\{p^2e,\,p^2b-b\}) \\
&=p^2e.
\end{align*}

     Suppose $m=p^2b-pk$ with
$p^2b-pk<\min\{p^2e,p^2b-b\}$.  We claim that
$p\nmid k$.  In fact if $p\mid k$ then
$k<b<\frac{pe}{p-1}$, so we have
\begin{equation} \nonumber
H\cap U_K^k=((K^{\times})^p\cap U_K^k)\cdot
(H\cap U_K^{k+1}).
\end{equation}
Since $pk\ge b+1$ and
$H\cap U_K^{k+1}\supset U_K^{b+1}$ it follows
from Corollary~\ref{sub} that
\begin{equation} \label{HUK}
\Lambda_p(H\cap U_K^k)=
\Lambda_p((K^{\times})^p\cap U_K^k)+
\Lambda_p(H\cap U_K^{k+1}).
\end{equation}
Since $p^2b-pk<p^2e$ we have $e+\frac{k}{p}\ge b+1$.
Therefore by Proposition~\ref{OK0} we see that
$\Lambda_p((K^{\times})^p\cap U_K^k)+\M_K^{b+1}$ is an
$\OO_{K_0}$-module.  Furthermore,
$\Lambda_p(H\cap U_K^{k+1})$ is a $\Zps$-module by the
definition of $k$.  Since $\Zps\subset\OO_{K_0}$ and
$\Lambda_p(H\cap U_K^{k+1})\supset\M_K^{b+1}$, it
follows from (\ref{HUK}) that
$\Lambda_p(H\cap U_K^k)$ is a $\Zps$-module.  This
contradicts the definition of $k$, so $p\nmid k$.

     Suppose $a_{pk}\in\M_K^{b+1}$.  Then for every
$\eta\in\mubold_{p^2-1}$ and $r\in\mubold_{q-1}$ we
have
\begin{equation} \label{etabeta}
\eta\beta_k(r)\equiv\beta_k(\eta r)\pmod{\pi_K^{b+1}}.
\end{equation}
If $\mubold_{p^2-1}\subset K$ this implies
$\eta\beta_k(r)\in\Lambda_p(H\cap U_K^k)$.  Since
$\Lambda_p(H\cap U_K^{k+1})$ is a $\Zps$-module it
follows that $\Lambda_p(H\cap U_K^k)$ is a
$\Zps$-module, contrary to assumption.  Therefore
$a_{pk}\not\in\M_K^{b+1}$ in this case.  If
$\mubold_{p^2-1}\not\subset K$ then $k=b$ and it follows
from (\ref{etabeta}) that the set
\begin{equation} \nonumber
S=\{r\in\mubold_{q-1}\cup\{0\}:
\beta_b(r)\equiv0\pmod{\M_K^{b+1}}\}
\end{equation}
is stable under multiplication by elements of
$\mubold_{p^2-1}$.  Hence $S=\{0\}$.  Since
\begin{equation} \nonumber
\beta_b(r+r')\equiv\beta_b(r)+\beta_b(r')\pmod{\M_K^{b+1}}
\end{equation}
for $r,r'\in\mubold_{q-1}\cup\{0\}$ this implies that
every coset of $\M_K^{b+1}$ in $\M_K^b$ is represented
by $\beta_b(r)$ for some $r\in\mubold_{q-1}\cup\{0\}$.
It follows that $\Lambda_p(H\cap U_K^b)=\M_K^b$, a
contradiction.  Hence $a_{pk}\not\in\M_K^{b+1}$ in this
case as well.

     Since $p\nmid k+p$, by (\ref{star}) we have
$\pi_Ka_{pk}=a_{p(k+p)}\in\M_K^{b+1}$.  Thus
$v_K(a_{pk})=b$.  Using (\ref{i1}) and (\ref{star}) we
get
\begin{align*}
i_1&=\min\left(\{p^2v_K(a_{pi})-pi:k\le i<k+p,\,p\nmid i\}
\cup\{p^2e,\,p^2b-b\}\right) \\
&=p^2b-pk.
\end{align*}
We conclude that $i_1=m$ in every case. \qed \medskip

\begin{remark}
Suppose $\mubold_{p^2-1}\subset K$.  Then it
follows from Corollary~\ref{OK} and class field theory
that all values of $k$ such that $b/p<k\le b$ and
$p\nmid k$ can be realized by extensions $L/K$
satisfying the conditions of Theorem~\ref{combining}.
\end{remark} \medskip

\begin{remark}
Using Theorem~\ref{combining} we obtain the bounds
$p^2b-pb\le i_1\le p^2b-b$.  These inequalities can
also be derived from Corollary~6.11 in \cite{heier}.
It follows from these bounds that the condition
$i_1>p^2b-pb$ is equivalent to $i_1\not=p^2b-pb$.
\end{remark}

\section{Kummer theory}

Let $p>2$ and let $K$ be a finite extension of $\Q_p$
which contains a primitive $p$th root of unity
$\zeta_p$.  Let $K^{ab}$ be a maximal abelian extension
of $K$ and let $L/K$ be a totally ramified
$(\Z/p\Z)^2$-subextension of $K^{ab}/K$ with a single
ramification break $b$.  In
\cite{necessity}, Byott and Elder gave a method
for computing the refined ramification break $b_*$ of
$L/K$ in terms of Kummer theory.  In this section we use
Vostokov's formula for the Kummer pairing to express
$b_*$ in terms of the index of inseparability $i_1$,
under the assumption that $i_1$ is not equal to
$p^2b-pb$.  The proof is based on a symmetry relation
involving the Kummer pairing and truncated
exponentiation.

     The Kummer pairing
$\langle\;,\,\rangle_p:K^{\times}\times K^{\times}
\ra\mubold_p$ is defined by
$\langle\alpha,\beta\rangle_p
=\sigma_{\beta}(\alpha^{1/p})/\alpha^{1/p}$, where
$\alpha^{1/p}\in K^{ab}$ is any
$p$th root of $\alpha$ and $\sigma_{\beta}$ is the
element of $\Gal(K^{ab}/K)$ that corresponds to $\beta$
under class field theory.  The Kummer pairing is
$\Z$-bilinear and skew-symmetric, with kernel
$(K^{\times})^p$ on the left and right (see for
instance Proposition~5.1 in \cite[IV]{FV}).
For $1\le i\le \frac{pe}{p-1}$ the
orthogonal complement of $U_K^i$ with respect to
$\langle\;,\,\rangle_p$ is $(U_K^i)^{\perp}
=(K^{\times})^p\cdot U_K^{\frac{pe}{p-1}-i+1}$ (see
\cite[\S1]{dalfurther}).

     Recall that $K_0/\Q_p$ is the maximum unramified
subextension of $K/\Q_p$.  In \cite{vos} Vostokov gave a
formula for computing $\langle\;,\,\rangle_p$ in terms
of residues of elements of
\begin{equation} \nonumber
K_0\{\!\{X\}\!\}=\left\{\sum_{n=-\infty}^{\infty}\,a_nX^n:
a_n\in K_0,\,\lim_{n\ra-\infty}\,v_{K_0}(a_n)=\infty,\,
\exists m\,\forall n\,v_{K_0}(a_n)\ge m\right\}.
\end{equation}
The set $K_0\{\!\{X\}\!\}$ has an obvious operation of
addition, and the conditions on the coefficients imply
that the natural multiplication on $K_0\{\!\{X\}\!\}$ is
also well-defined.  These operations make
$K_0\{\!\{X\}\!\}$ a field.  Let
$\OO_{K_0}\{\!\{X\}\!\}$ denote the subring of
$K_0\{\!\{X\}\!\}$ consisting of series whose
coefficients lie in $\OO_{K_0}$.  Also let
$\Res(\psi(X))$ denote the coefficient of $X^{-1}$ in
$\psi(X)\in K_0\{\!\{X\}\!\}$.

     For each $\alpha\in U_K^1$ choose
$\tilde{\alpha}(X)\in\OO_{K_0}[[X]]$ so that
$\tilde{\alpha}(0)=1$ and
$\tilde{\alpha}(\pi_K)=\alpha$.  Of course there
are many series $\tilde{\alpha}(X)$ with this property,
but for our purposes it will not matter which we choose.
Let $\phi:K_0\ra K_0$ be the $p$-Frobenius map and
define
$\tilde{\alpha}^{\Delta}(X)=\tilde{\alpha}^{\phi}(X^p)$
and $l(\tilde{\alpha})=\log(\tilde{\alpha})-
p^{-1}\log(\tilde{\alpha}^{\Delta})$,
where
\begin{equation} \nonumber
\log(1+\psi(X))=\tst\psi(X)-\frac12\psi(X)^2+
\frac13\psi(X)^3-\dots
\end{equation}
for $\psi(X)\in XK_0[[X]]$.  By Proposition~2.2 in
\cite[VI]{FV} we have
$l(\tilde{\alpha})\in X\OO_{K_0}[[X]]$.

     Let $\alpha,\beta\in U_K^1$.
Following \cite[p.\,241]{FV} we define
\begin{alignat*}{2}
\Phi_{\alpha,\beta}(X)
&=\frac{\tilde{\alpha}'}{\tilde{\alpha}}\cdot
l(\tilde{\beta})-
\frac{(\tilde{\beta}^{\Delta})'}{p\tilde{\beta}^{\Delta}}
\cdot l(\tilde{\alpha}).
\end{alignat*}
Then $\Phi_{\alpha,\beta}(X)\in\OO_{K_0}[[X]]$.  Let
$s(X)=\tilde{\zeta}_p(X)^p-1$.  Then by Proposition~3.1 in
\cite[VI]{FV}, $s(X)$ is a unit in
$\OO_{K_0}\{\!\{X\}\!\}$.  Since $p>2$ and
$\alpha,\beta\in U_K^1$, by Theorem~4 in \cite[VII]{FV}
we have
\begin{equation} \label{vost}
\langle\alpha,\beta\rangle_p
=\zeta_p^{\Tr_{K_0/\Q_p}(\Res(\Phi_{\alpha,\beta}/s))}.
\end{equation}

\begin{theorem} \label{sym}
Let $p>2$ and let $K$ be a finite extension of $\,\Q_p$
which contains a primitive $p$th root of unity.  Let
$i,j$ be positive integers such that
$i+pj>\frac{pe}{p-1}$ and $pi+j>\frac{pe}{p-1}$.  Let
$\alpha\in U_K^i$, $\beta\in U_K^j$, and
$\eta\in\OO_{K_0}$.  Then
$\langle\alpha^{[\eta]},\beta\rangle_p=
\langle\alpha,\beta^{[\eta]}\rangle_p$.
\end{theorem}

\proof By the linearity and continuity of the Kummer
pairing we may assume that $\alpha=E_p(u\pi_K^c)$,
$\beta=E_p(v\pi_K^d)$, $\tilde{\alpha}(X)=E_p(uX^c)$,
and $\tilde{\beta}(X)=E_p(vX^d)$ with
$u,v\in{\mubold}_{q-1}$, $c\ge i$, and $d\ge j$.  It
follows from (\ref{Ep}) that $l(\tilde{\alpha}(X))=uX^c$
and $l(\tilde{\beta}(X))=vX^d$.  Using (\ref{Ep}) and
Lemma~\ref{alpha} we get
\begin{alignat*}{2}
\frac{\tilde{\alpha}'(X)}{\tilde{\alpha}(X)}
&\equiv cuX^{c-1}&&\pmod{X^{pc-1}} \\
\frac{(\tilde{\beta}^{\Delta})'(X)}{p\tilde{\beta}^{\Delta}(X)}
&\equiv0&&\pmod{X^{pd-1}} \\
\frac{(\tilde{\alpha}(X)^{[\eta]})'}{\tilde{\alpha}(X)^{[\eta]}}
&\equiv c(\eta u)X^{c-1}&&\pmod{X^{pc-1}} \\
l(\tilde{\beta}(X)^{[\eta]})&\equiv\eta vX^d
&&\pmod{X^{pd}}.
\end{alignat*}

     Note that $\tilde{\alpha}(X)^{[\eta]}$,
$\tilde{\beta}(X)^{[\eta]}$ are elements of
$1+X\OO_{K_0}[[X]]$ such that
$\tilde{\alpha}(\pi_K)^{[\eta]}=\alpha^{[\eta]}$,
$\tilde{\beta}(\pi_K)^{[\eta]}=\beta^{[\eta]}$.  Hence
we may take $\widetilde{\alpha^{[\eta]}}(X)
=\tilde{\alpha}(X)^{[\eta]}$ and
$\widetilde{\beta^{[\eta]}}(X)=\tilde{\beta}(X)^{[\eta]}$.
Using the computations from the preceding paragraph and
the lower bounds for $i+pj$ and $pi+j$ we get
\begin{alignat}{2}
\Phi_{\alpha,\beta}(X)&\equiv
\frac{\tilde{\alpha}'}{\tilde{\alpha}}\cdot
l(\tilde{\beta})&&\pmod{X^{\frac{pe}{p-1}}} \nonumber \\
\Phi_{\alpha^{[\eta]},\beta}(X)
&\equiv c(\eta u)vX^{c+d-1}&&\pmod{X^{\frac{pe}{p-1}}}
\label{Phiaomegab} \\
\Phi_{\alpha,\beta^{[\eta]}}(X)
&\equiv cu(\eta v)X^{c+d-1}&&\pmod{X^{\frac{pe}{p-1}}}.
\label{Phiabomega}
\end{alignat}
It follows from Proposition~3.1 in \cite[VI]{FV} that
the image of $s(X)\in\OO_{K_0}\{\!\{X\}\!\}^{\times}$ in
\begin{equation} \nonumber
(\OO_{K_0}/\M_{K_0})((X))\cong\Fq((X))
\end{equation}
has $X$-valuation $\frac{pe}{p-1}$.  Therefore by
(\ref{Phiaomegab}) and (\ref{Phiabomega}) we have
\begin{equation} \nonumber
\frac{\Phi_{\alpha^{[\eta]},\beta}(X)
-\Phi_{\alpha,\beta^{[\eta]}}(X)}{s(X)}
=\gamma(X)+p\delta(X)
\end{equation}
for some $\gamma(X)\in\OO_{K_0}[[X]]$ and
$\delta(X)\in\OO_{K_0}\{\!\{X\}\!\}$.  It follows that
\begin{alignat*}{2}
\Res\left(\frac{\Phi_{\alpha^{[\eta]},\beta}(X)}{s(X)}\right)
&\equiv
\Res\left(\frac{\Phi_{\alpha,\beta^{[\eta]}}(X)}{s(X)}\right)
&&\pmod{\M_{K_0}}.
\end{alignat*}
Therefore by (\ref{vost}) we get
$\langle\alpha^{[\eta]},\beta\rangle_p=
\langle\alpha,\beta^{[\eta]}\rangle_p$. \qed

\begin{cor} \label{perp}
Let $K$, $i$, $j$ satisfy the hypotheses of
Theorem~\ref{sym}.  Let
$A$ be a subgroup of $U_K^i$ such that $A$ contains
$U_K^{pi}$ and $\Lambda_p(A)$ is a $\Zps$-module.  Then
$\Lambda_p(A^{\perp}\cap U_K^j)$ is a $\Zps$-module.
\end{cor}

\proof Let $\alpha\in A$.  By Corollary~\ref{module}
we have $\alpha^{[\eta]}\in A$ for every
$\eta\in\mubold_{p^2-1}$.  Hence for
$\beta\in A^{\perp}\cap U_K^j$ we see that
$\langle\alpha,\beta^{[\eta]}\rangle_p
=\langle\alpha^{[\eta]},\beta\rangle_p=1$.  Since this
holds for every $\alpha\in A$ we get
$\beta^{[\eta]}\in A^{\perp}\cap U_K^j$.  Since
$pj\ge\frac{pe}{p-1}-i+1$ we have
$A^{\perp}\cap U_K^j\supset U_K^{pj}$.  Therefore
it follows from Corollary~\ref{module} that
$\Lambda_p(A^{\perp}\cap U_K^j)$ is a $\Zps$-module.
\qed \medskip

     Recall that $H=\N_{L/K}(L^{\times})$ is the
subgroup of $K^{\times}$ that corresponds to $L/K$ under
class field theory, and let
$R=(L^{\times})^p\cap K^{\times}$ denote the subgroup of
$K^{\times}$ that corresponds to $L/K$ under Kummer
theory.  Then $R$ contains $(K^{\times})^p$, and it
follows from the basic properties of the Kummer pairing
that $R=H^{\perp}$ and $H=R^{\perp}$.  Furthermore,
$R/(K^{\times})^p$ and $K^{\times}/H$ are both
elementary abelian $p$-groups of rank 2.  Let
$R_0=R\cap U_K^{\frac{pe}{p-1}-b}$.  Since the only
ramification break of $L/K$ is $b$ we see that
$R=R_0\cdot(K^{\times})^p$ and
\begin{equation} \nonumber
R_0/((K^{\times})^p\cap U_K^{\frac{pe}{p-1}-b})
\cong R/(K^{\times})^p
\end{equation}
(cf.\ \cite{dalfurther}).

     For $a\in\OO_K$ we let
$\overline{a}=a+\M_K^{\frac{pe}{p-1}-b+1}$ denote the
image of $a$ in $\OO_K/\M_K^{\frac{pe}{p-1}-b+1}$.
Then $\Rbar_0\cong R/(K^{\times})^p$ is an
elementary abelian $p$-group of rank 2.  Let
$1+\rho_1,1+\rho_2$ be elements of $R_0$ such that
$\overline{1+\rho_1}$, $\overline{1+\rho_2}$ generate
$\Rbar_0$.  Then
$v_K(\rho_1)=v_K(\rho_2)=\frac{pe}{p-1}-b$.  Let
$\theta\in\mubold_{q-1}$ be such that
$\theta\equiv\rho_2/\rho_1\pmod{\M_K}$.  Then
$\theta\not\in\mubold_{p-1}$ and
\begin{alignat*}{2}
(1+\rho_1)^{[\theta]}
&\equiv1+\rho_2&&\pmod{\M_K^{\frac{pe}{p-1}-b+1}}.
\end{alignat*}
Let $s\le\frac{pe}{p-1}$ be maximum such that
$(1+\rho_1)^{[\theta]}\in R_0\cdot U_K^s$, and
set $t=\frac{pe}{p-1}-s$.  Then by
\cite[Prop.\,10]{necessity} we have
\begin{equation} \label{bstar}
b_*=pb-\max\{pt-b,\,(p^2-1)b-p^2e,\,0\}.
\end{equation}

\begin{lemma} \label{Zpsmodule}
Let $p>2$ and assume that $K$ contains a primitive $p$th
root of unity.  Let $L/K$ be a totally ramified
$(\Z/p\Z)^2$-subextension of $K^{ab}/K$ with a single
ramification break $b$.  Then the following are equivalent:
\begin{enumerate}
\item $\theta\in\mubold_{p^2-1}$.
\item $\Lambda_p(R_0)+\M_K^{\frac{pe}{p-1}-b+1}$ is a
$\Zps$-module.
\item $\Lambda_p(H\cap U_K^b)$ is a $\Zps$-module.
\item $i_1>p^2b-pb$.
\end{enumerate}
\end{lemma}

\proof To prove the equivalence of the first two
statements we note that $\overline{\Lambda_p(1+\rho_1)}$
and $\overline{\Lambda_p(1+\rho_2)}
=\theta\cdot\overline{\Lambda_p(1+\rho_1)}$ generate the
rank-2 elementary abelian $p$-group
$\overline{\Lambda_p(R_0)}$.  Hence
$\theta$ lies in $\mubold_{p^2-1}$ if and only if
$\overline{\Lambda_p(R_0)}$ is a vector space over
$\Fps$, which holds if and only if
$\Lambda_p(R_0)+\M_K^{\frac{pe}{p-1}-b+1}$ is a
$\Zps$-module.  The equivalence of statements 3 and 4
follows from Theorem~\ref{combining}.  To prove the
equivalence of statements 2 and 3 we observe that if
$\Lambda_p(R_0)+\M_K^{\frac{pe}{p-1}-b+1}$ is a
$\Zps$-module then it follows from Corollary~\ref{perp}
that
\begin{equation} \nonumber
\Lambda_p((R_0\cdot U_K^{\frac{pe}{p-1}-b+1})^{\perp}
\cap U_K^b)=\Lambda_p(H\cap U_K^b)
\end{equation}
is a $\Zps$-module.  Conversely, if
$\Lambda_p(H\cap U_K^b)$ is a $\Zps$-module then it
follows from Corollary~\ref{perp} that
\begin{align*}
\Lambda_p((H\cap U_K^b)^{\perp}\cap
U_K^{\frac{pe}{p-1}-b})&=\Lambda_p(R_0\cdot
U_K^{\frac{pe}{p-1}-b+1}) \\
&=\Lambda_p(R_0)+\M_K^{\frac{pe}{p-1}-b+1}
\end{align*}
is a $\Zps$-module. \qed \medskip

     For the rest of this paper we restrict our
attention to extensions $L/K$ which satisfy the
conditions of Lemma~\ref{Zpsmodule}.  Our goal is to
compute $b_*$ in terms of $i_1$ for this class of
extensions.  The following proposition will allow us to
make a connection between $\Lambda_p(R_0)$ and the
definition of $s$.

\begin{prop} \label{Zps}
Let $L/K$ be an extension which satisfies the conditions
of Lemma~\ref{Zpsmodule}, and let $i$ satisfy
$1\le i\le p(\frac{pe}{p-1}-b)$ and
$i\le\frac{pe}{p-1}-\lfloor\frac{b}{p}\rfloor$.  Then
$(1+\rho_1)^{[\theta]}\in R_0\cdot U_K^i$ if and only if
$\Lambda_p(R_0)+\M_K^i$ is a $\Zps$-module.
\end{prop}

\proof If $i\le\frac{pe}{p-1}-b$ then both statements
are certainly true, so we assume $i>\frac{pe}{p-1}-b$.
If $\Lambda_p(R_0)+\M_K^i$ is a $\Zps$-module
then it follows from Proposition~\ref{hom} that
$(1+\rho_1)^{[\theta]}\in R_0\cdot U_K^i$.  Conversely,
suppose that $(1+\rho_1)^{[\theta]}\in R_0\cdot U_K^i$.
Thanks to the upper bounds on $i$, the
hypotheses of Proposition~\ref{OK0} are satisfied with
$j=\frac{pe}{p-1}-b$.  It follows that
$\Lambda_p((K^{\times})^p\cap U_K^{\frac{pe}{p-1}-b})
+\M_K^i$ is an $\OO_{K_0}$-module, and hence a
$\Zps$-module.  By Proposition~\ref{hom} we have
$\theta\cdot\Lambda_p(1+\rho_1)\in\Lambda_p(R_0)+\M_K^i$.
Therefore the rank-2 elementary abelian $p$-group
\begin{equation} \label{quot2}
(\Lambda_p(R_0)+\M_K^i)/
(\Lambda_p((K^{\times})^p\cap U_K^{\frac{pe}{p-1}-b})
+\M_K^i)
\end{equation}
is generated by the cosets represented by
$\Lambda_p(1+\rho_1)$ and $\theta\cdot\Lambda_p(1+\rho_1)$.
Since
$\theta\in\mubold_{p^2-1}\smallsetminus\mubold_{p-1}$,
it follows that (\ref{quot2}) is a vector space over
$\Fps$.  We conclude that $\Lambda_p(R_0)+\M_K^i$
is a $\Zps$-module. \qed \medskip

     We now reformulate the Byott-Elder formula for
$b_*$ in terms of $\Lambda_p(R_0)$.

\begin{theorem} \label{sprime}
Let $L/K$ be an extension which satisfies the conditions
of Lemma~\ref{Zpsmodule}, let $R$ be the subgroup of
$K^{\times}$ that corresponds to $L/K$ under Kummer
theory, and set $R_0=R\cap U_K^{\frac{pe}{p-1}-b}$.  Let
$s'\le\frac{pe}{p-1}$ be maximum such that
$\Lambda_p(R_0)+\M_K^{s'}$ is a $\Zps$-module and set
$t'=\frac{pe}{p-1}-s'$.  Then
\begin{equation}
b_*=pb-\max\{pt'-b,\,(p^2-1)b-p^2e,\,0\}.
\end{equation}
\end{theorem}

\proof Recall that $t=\frac{pe}{p-1}-s$, where $s$ is
the smallest nonnegative integer such that
$(1+\rho_1)^{[\theta]}\in R_0\cdot U_K^s$.  Set
\begin{align*}
M&=\max\{pt-b,\,(p^2-1)b-p^2e,\,0\} \\
M'&=\max\{pt'-b,\,(p^2-1)b-p^2e,\,0\}.
\end{align*}
By (\ref{bstar}) we have $b_*=pb-M$.  Therefore to prove
the theorem it suffices to show that $M'=M$.  We
divide the proof into three cases, depending on the
value of $M$.

     If $M=(p^2-1)b-p^2e$ then $t\le p(b-e)$, and hence
$(1+\rho_1)^{[\theta]}\in R_0\cdot
U_K^{\frac{pe}{p-1}-p(b-e)}$.  Since
$(p^2-1)b-p^2e\ge0$ we have
\begin{align*}
p\left(\frac{pe}{p-1}-b\right)=
\frac{pe}{p-1}-p(b-e)
&\le\frac{pe}{p-1}-\left\lfloor\frac{b}{p}\right\rfloor.
\end{align*}
Therefore by Proposition~\ref{Zps} we see that
$\Lambda_p(R_0)+\M_K^{\frac{pe}{p-1}-p(b-e)}$ is a
$\Zps$-module.  Hence $t'\le p(b-e)$, so $M'=M$ in
this case.  

     If $M=0$ then $t\le\lfloor\frac{b}{p}\rfloor$ and
hence $(1+\rho_1)^{[\theta]}\in R_0\cdot
U_K^{\frac{pe}{p-1}-\lfloor\frac{b}{p}\rfloor}$.
Since $(p^2-1)b-p^2e\le0$ we have
$p(\frac{pe}{p-1}-b)\ge\frac{pe}{p-1}-\lfloor
\frac{b}{p}\rfloor$.  Therefore by
Proposition~\ref{Zps} we see
that $\Lambda_p(R_0)+\M_K^{\frac{pe}{p-1}-\lfloor
\frac{b}{p}\rfloor}$ is a $\Zps$-module.  Hence
$t'\le\lfloor \frac{b}{p}\rfloor$, so $pt'\le b$.
It follows that $M'=M$ in this case.

     If $M=pt-b>\max\{(p^2-1)b-p^2e,0\}$ then
$t>p(b-e)$ and $t>\frac{b}{p}$.  Hence
$s<p(\frac{pe}{p-1}-b)$ and
$s<\frac{pe}{p-1}-\lfloor\frac{b}{p}\rfloor$.
Since $(1+\rho_1)^{[\theta]}\in R_0\cdot U_K^s$
and $(1+\rho_1)^{[\theta]}\not\in R_0\cdot U_K^{s+1}$,
it follows from Proposition~\ref{Zps} that
$\Lambda_p(R_0)+\M_K^s$ is a $\Zps$-module, but
$\Lambda_p(R_0)+\M_K^{s+1}$ is not.  Therefore
$s'=s$, so $M'=M$ in this case as well. \qed \medskip

     Now that we have formulas for computing $b_*$ and
$i_1$ in terms of $\Lambda_p(R_0)$, we can determine the
relationship between these two invariants.

\begin{theorem}
Let $p>2$ and let $K$ be a finite extension of $\Q_p$
which contains a primitive $p$th root of unity.  Let
$L/K$ be a
totally ramified $(\Z/p\Z)^2$-extension with a single
ramification break $b$.  Assume that the index of
inseparability $i_1$ of $L/K$ is not equal to $p^2b-pb$.
Then the refined ramification break
$b_*$ of $L/K$ is given by $b_*=i_1-p^2b+pb+b$.
\end{theorem}

\proof As above we let $H$ denote the subgroup of
$K^{\times}$ that corresponds to the extension $L/K$
under class field theory.  By Theorem~\ref{combining} we
have
\begin{equation} \label{min}
i_1=\min\{p^2b-pk,\,p^2e,\,p^2b-b\},
\end{equation}
where $k$ is the smallest nonnegative integer such
that $\Lambda_p(H\cap U_K^{k+1})$ is a $\Zps$-module.
Let $R$ be the subgroup of $K^{\times}$ that
corresponds to $L/K$ under Kummer theory and set
$R_0=R\cap U_K^{\frac{pe}{p-1}-b}$.  Recall that
$R$ is equal to the orthogonal complement $H^{\perp}$ of
$H$ with respect to the Kummer pairing
$\langle\;,\,\rangle_p$.
In addition, since $R=R_0\cdot(K^{\times})^p$ we have
$R_0^{\perp}=R^{\perp}=H$.  As in Theorem~\ref{sprime}
we let $t'$ be the smallest nonnegative integer such
that $\Lambda_p(R_0)+\M_K^{\frac{pe}{p-1}-t'}$ is a
$\Zps$-module.

     Suppose $i_1=p^2b-b$.  Then
\begin{align*}
\Lambda_p((H\cap
U_K^{\lfloor\frac{b}{p}\rfloor+1})^{\perp}\cap
U_K^{\frac{pe}{p-1}-b})&=
\Lambda_p((R\cdot
U_K^{\frac{pe}{p-1}-\lfloor\frac{b}{p}\rfloor})\cap
U_K^{\frac{pe}{p-1}-b}) \\
&=\Lambda_p(R_0\cdot
U_K^{\frac{pe}{p-1}-\lfloor\frac{b}{p}\rfloor}).
\end{align*}
Since $p(\frac{pe}{p-1}-b)\ge\frac{pe}{p-1}-
\lfloor\frac{b}{p}\rfloor$, it follows from
Corollary~\ref{sub} that
\begin{equation} \label{pR0}
\Lambda_p((H\cap
U_K^{\lfloor\frac{b}{p}\rfloor+1})^{\perp}\cap
U_K^{\frac{pe}{p-1}-b})=\Lambda_p(R_0)+
\M_K^{\frac{pe}{p-1}-\lfloor\frac{b}{p}\rfloor}.
\end{equation}
Since $\lfloor\frac{b}{p}\rfloor+1>\frac{b}{p}\ge
p(b-e)$, we have
\begin{align*}
\left(\left\lfloor\frac{b}{p}\right\rfloor+1\right)
+p\left(\frac{pe}{p-1}-b\right)&>\frac{pe}{p-1} \\
p\left(\left\lfloor\frac{b}{p}\right\rfloor+1\right)
+\left(\frac{pe}{p-1}-b\right)&>\frac{pe}{p-1}.
\end{align*}
Therefore by (\ref{pR0}) and Corollary~\ref{perp} with
$A=H\cap U_K^{\lfloor\frac{b}{p}\rfloor+1}$,
$i=\lfloor\frac{b}{p}\rfloor+1$, and
$j=\frac{pe}{p-1}-b$ we see that $\Lambda_p(R_0)+
\M_K^{\frac{pe}{p-1}-\lfloor\frac{b}{p}\rfloor}$
is a $\Zps$-module.  Hence $t'\le
\lfloor\frac{b}{p}\rfloor$.  Since $(p^2-1)b-p^2e\le0$,
it follows from Theorem~\ref{sprime} that $b_*=pb$ in
this case.

     Suppose $i_1=p^2e$.  Then
\begin{align*}
\Lambda_p((H\cap U_K^{p(b-e)+1})^{\perp}\cap
U_K^{\frac{pe}{p-1}-b})&=
\Lambda_p((R\cdot U_K^{\frac{pe}{p-1}-p(b-e)})\cap
U_K^{\frac{pe}{p-1}-b}) \\
&=\Lambda_p(R_0\cdot U_K^{\frac{pe}{p-1}-p(b-e)}).
\end{align*}
Since $b>p(b-e)$ and
$p(\frac{pe}{p-1}-b)=\frac{pe}{p-1}-p(b-e)$ it follows
from Corollary~\ref{sub} that
\begin{equation} \label{pR1}
\Lambda_p((H\cap U_K^{p(b-e)+1})^{\perp}\cap
U_K^{\frac{pe}{p-1}-b})
=\Lambda_p(R_0)+\M_K^{\frac{pe}{p-1}-p(b-e)}.
\end{equation}
Since $p^2b-b\ge p^2e$ we have
\begin{align*}
(p(b-e)+1)+p\left(\frac{pe}{p-1}-b\right)&>\frac{pe}{p-1} \\
p(p(b-e)+1)+\left(\frac{pe}{p-1}-b\right)&>\frac{pe}{p-1}.
\end{align*}
Therefore it follows from (\ref{pR1}) and
Corollary~\ref{perp} with $A=H\cap U_K^{p(b-e)+1}$,
$i=p(b-e)+1$, and $j=\frac{pe}{p-1}-b$ that
$\Lambda_p(R_0)+\M_K^{\frac{pe}{p-1}-p(b-e)}$
is a $\Zps$-module.  Hence $t'\le p(b-e)$.  Since
$(p^2-1)b-p^2e\ge0$, it follows from Theorem~\ref{sprime}
that $b_*=p^2(e-b)+pb+b$ in this case.

     Suppose $i_1=p^2b-pk<\min\{p^2b-b,p^2e\}$.
Since $H\supset U_K^{b+1}$ we have $k\le b$, so
$R_0\cdot U_K^{\frac{pe}{p-1}-k}$ is contained in
$U_K^{\frac{pe}{p-1}-b}$.  Hence
\begin{align*}
\Lambda_p((H\cap U_K^{k+1})^{\perp}\cap
U_K^{\frac{pe}{p-1}-b})&=
\Lambda_p((R\cdot U_K^{\frac{pe}{p-1}-k})\cap
U_K^{\frac{pe}{p-1}-b}) \\
&=\Lambda_p(R_0\cdot U_K^{\frac{pe}{p-1}-k}).
\end{align*}
Since $k>p(b-e)$ we have
$p(\frac{pe}{p-1}-b)>\frac{pe}{p-1}-k$.  Therefore by
Corollary~\ref{sub} we get
\begin{equation} \label{pR2}
\Lambda_p((H\cap U_K^{k+1})^{\perp}\cap
U_K^{\frac{pe}{p-1}-b})
=\Lambda_p(R_0)+\M_K^{\frac{pe}{p-1}-k}.
\end{equation}
It follows from the
inequalities $k>p(b-e)$ and $pk>b$ that
\begin{align*}
k+p\left(\frac{pe}{p-1}-b\right)&>\frac{pe}{p-1} \\
pk+\left(\frac{pe}{p-1}-b\right)&>\frac{pe}{p-1}.
\end{align*}
Therefore by (\ref{pR2}) and Corollary~\ref{perp} with
$A=H\cap U_K^{k+1}$,
$i=k+1$, and $j=\frac{pe}{p-1}-b$ we see that
$\Lambda_p(R_0)+\M_K^{\frac{pe}{p-1}-k}$
is a $\Zps$-module.

     Suppose that
$\Lambda_p(R_0)+\M_K^{\frac{pe}{p-1}-k+1}$ is also a
$\Zps$-module.  Then by Corollary~\ref{perp} with
$A=R_0\cdot U_K^{\frac{pe}{p-1}-k+1}$,
$i=\frac{pe}{p-1}-b$, and $j=k$ we see that
\begin{align*}
\Lambda_p((R_0\cdot U_K^{\frac{pe}{p-1}-k+1})^{\perp}\cap U_K^k)
&=\Lambda_p(H\cap(K^{\times})^pU_K^k\cap U_K^k) \\
&=\Lambda_p(H\cap U_K^k)
\end{align*}
is a $\Zps$-module.  Since $k\ge1$ this contradicts the
definition of $k$.  Hence $\Lambda_p(R_0\cdot
U_K^{\frac{pe}{p-1}-k+1})$ is not a $\Zps$-module, so
$t'=k$.  Since $pk-b>\max\{(p^2-1)b-p^2e,0\}$ we
get $b_*=pb-pk+b$ by Theorem~\ref{sprime}.  By comparing
our formulas for $b_*$ with (\ref{min}) we find that
$b_*=i_1-p^2b+pb+b$ in all three cases. \qed \medskip

\begin{remark}
If $i_1=p^2b-pb$ then $b_*$ can take any of the values
allowed by Theorem~5 in \cite{necessity}.  On the other
hand, for a given $b_*$ we have either $i_1=p^2b-pb$ or
$i_1=b_*+p^2b-pb-b$.
\end{remark}

\end{document}